\documentclass[reqno,11pt]{amsart}

\usepackage{xcolor}

\usepackage{amssymb}
\usepackage[all,cmtip]{xy}
\usepackage{comment}
\newtheorem{Proposition}{Proposition}[section]

\newtheorem{Definition}[Proposition]{Definition}
\newtheorem{Lemma}[Proposition]{Lemma}
\newtheorem{Theorem}[Proposition]{Theorem}

\DeclareMathOperator{\Val}{Val}
\DeclareMathOperator{\Area}{Area}

\DeclareMathOperator{\Curv}{Curv}
\DeclareMathOperator{\Ang}{Ang}
\DeclareMathOperator{\glob}{glob}
\DeclareMathOperator{\PD}{PD}
\DeclareMathOperator{\pd}{pd}
\newcommand{\R}{\mathbb{R}}
\newcommand{\C}{\mathbb{C}}
\DeclareMathOperator{\CP}{\C P}
\title{Dual curvature measures in hermitian integral geometry}
\author{Andreas Bernig, Joseph H.G. Fu, and Gil Solanes}
 
\email{bernig@math.uni-frankfurt.de}
\email{fu@math.uga.edu}
\email{solanes@mat.uab.cat}

\address{Institut f\"ur Mathematik, Goethe-Universit\"at Frankfurt,
Robert-Mayer-Str. 10, 60054 Frankfurt, Germany}
\address{ Department of Mathematics, 
University of Georgia, 
Athens, GA 30602, USA}
\address{Departament de Matem\`atiques, Universitat Aut\`onoma de Barcelona, 08193 Bellaterra, Spain}

\thanks{A.B. was supported by DFG grant BE 2484/5-2. J.H.G.F. was supported by NSF grant DMS-1406252. G.S. is a  Serra H\'unter Fellow and was supported by FEDER-MINECO grant MTM2015-66165-P}

\begin{document}
\begin{abstract}
The local kinematic formulas on complex space forms induce the structure of a commutative algebra on the space $\mathrm{Curv}^{\mathrm{U}(n)*}$ of dual unitarily invariant curvature measures. Building on the recent results from integral geometry in complex space forms, we describe this algebra structure explicitly as a polynomial algebra. This is a short way to encode all local kinematic formulas. We then characterize the invariant valuations on complex space forms leaving the space of invariant angular curvature measures fixed. 
 \end{abstract}

\maketitle 
%---------------------------------------------------------------------

\section{Introduction}

Let $\CP_\lambda^n$ denote the complex space form of holomorphic sectional curvature $4\lambda$ and $G_\lambda$  its holomorphic isometry group. By $\mathcal{C}(\CP^n_\lambda)^{G_\lambda}$ we denote the space of $G_\lambda$-invariant smooth curvature measures on $\CP^n_\lambda$ (see below for the definition). 

The space $\mathcal{C}(\CP^n_\lambda)^{G_\lambda}$ is finite-dimensional,  and several geometrically meaningful bases were used in \cite{bernig_fu_solanes}. Let $\Phi_1,\ldots,\Phi_m$ be  such a basis. Then there are local kinematic formulas  (cf. \cite{fu90} or \cite{fu_pokorny_rataj})
\begin{displaymath}
 \int_{G_\lambda} \Phi_j(P_1 \cap gP_2, \beta_1 \cap g\beta_2) dg = \sum_{k,l} c_{k,l}^j \Phi_k(P_1,\beta_1) \Phi_l(P_2,\beta_2).
\end{displaymath}
Here $P_1,P_2$ are compact submanifolds with corners and $\beta_1,\beta_2$ are Borel subsets of $\CP^n_\lambda$. Finding the constants $c_{k,l}^j$ is a non-trivial question which could be solved only recently \cite{bernig_fu_hig,bernig_fu_solanes}.

By the transfer principle \cite{bernig_fu_solanes}, the spaces $\mathcal{C}(\CP^n_\lambda)^{G_\lambda}$ are naturally identified with the space $\Curv^{\mathrm{U}(n)}$ of smooth, translation  and $\mathrm U(n)$-invariant curvature measures on the hermitian space $\C^n$ and under this identification, the local kinematic formulas are independent of $\lambda$. We may therefore define an operator
\begin{displaymath}
 K: \Curv^{\mathrm{U}(n)} \to \Curv^{\mathrm{U}(n)} \otimes \Curv^{\mathrm{U}(n)},\quad \Phi_j \mapsto \sum_{k,l} c_{k,l}^j \Phi_k \otimes \Phi_l,
\end{displaymath}
which makes $\Curv^{\mathrm{U}(n)}$ into a cocommutative, coassociative coalgebra. 

Stated otherwise, the dual space $\Curv^{\mathrm{U}(n)*}$ becomes a commutative associative algebra with respect to the product 
\begin{displaymath}
 K^*:\Curv^{\mathrm{U}(n)*} \otimes \Curv^{\mathrm{U}(n)*} \to \Curv^{\mathrm{U}(n)*}.
\end{displaymath}

The knowledge of the algebra structure on $\Curv^{\mathrm{U}(n)*}$ is equivalent to the knowledge of the local kinematic formulas. Hence the description of this structure is a short and elegant way of stating the local kinematic formulas. This will be achieved by the first main theorem of this paper. 

Before stating it, let us introduce some notation and recall a result by Fu \cite{fu06}. Let $\Val^{\mathrm{U}(n)}$ be the space of continuous, translation invariant and $\mathrm{U}(n)$-invariant valuations (see Section \ref{sec_background} for the definition of valuations). Let $\phi_1,\ldots,\phi_r$ be a basis of $\Val^{\mathrm{U}(n)}$. Then there are global kinematic formulas
\begin{displaymath}
 \int_{\overline{\mathrm{U}(n)}} \phi_j( P_1 \cap \bar g P_2) d\bar g=\sum_{k,l} \tilde c_{k,l}^j \phi_k(P_1)\phi_l(P_2).
\end{displaymath}
The operator $k:\Val^{\mathrm{U}(n)} \to \Val^{\mathrm{U}(n)} \otimes \Val^{\mathrm{U}(n)}$ defined by $\phi_j \mapsto \sum_{k,l} \tilde c_{k,l}^j \phi_k \otimes \phi_l$ induces the structure of a commutative, associative algebra 
\begin{displaymath}
k^*: \Val^{\mathrm{U}(n)^*} \otimes \Val^{\mathrm{U}(n)^*} \to \Val^{\mathrm{U}(n)^*}.
\end{displaymath}
In this case, we may identify $\Val^{\mathrm{U}(n)}$ and $\Val^{\mathrm{U}(n)^*}$ by the Alesker-Poincar\'e pairing and the resulting algebra structure on $\Val^{\mathrm{U}(n)}$ is the Alesker product of valuations defined in \cite{alesker04_product}. 

\begin{Theorem}[\cite{fu06}]
 The algebra $\Val^{\mathrm{U}(n)}$ is given by 
\begin{displaymath}
 \Val^{\mathrm{U}(n)} \cong \C[t,s]/(f_{n+1},f_{n+2}), 
\end{displaymath}
where $t,s$ are variables of degrees $1$ and $2$ respectively, and $f_k$ is the part of total degree $k$ in the series expansion of $\log(1+t+s)$.
\end{Theorem}

Let $t,s,v$ be variables of degrees $1$, $2$ and $3$ respectively and $u:=4s-t^2$. Wannerer \cite{wannerer_unitary_module} defined the polynomial $q_n$ as the $n$-homogeneous part in the expansion of $-\frac{1}{(1+t+s)^2}$. These polynomials appear in the description of the algebra $\Area^{\mathrm U(n)*}$ of unitarily invariant dual area measures \cite{wannerer_area_measures,wannerer_unitary_module}.

Our first main theorem gives an algebraic way of encoding all local kinematic formulas, and, moreover, answers \cite[Question  7.3]{bernig_fu_solanes}. 

\begin{Theorem} \label{mainthm_algebra}
There is an algebra isomorphism
\begin{displaymath}
\Curv^{U(n)*} \cong \C[t,s,v]/(f_{n+1},f_{n+2},q_{n-1}v,q_nv, (v+tu)^2). 
\end{displaymath}
\end{Theorem}

Let us now describe our second main theorem. Let $M$ be a Riemannian manifold. The space $\mathcal{V}(M)$ of smooth valuations on $M$ admits an algebra structure, and the space $\mathcal{C}(M)$ of smooth curvature measures on $M$ is a module over  this algebra \cite{bernig_fu_solanes}. In \cite{bernig_fu_solanes} we defined the subspace $\mathcal{A}(M) \subset \mathcal{C}(M)$ of {\it angular curvature measures}. A smooth valuation $\mu$ on $M$ with the property that $\mu \cdot \mathcal{A}(M) \subset \mathcal{A}(M)$ is called {\it angular valuation}; the corresponding space is denoted by $\mathfrak{a}(M) \subset \mathcal{V}(M)$. The Lipschitz-Killing algebra of $M$ (\cite[Subsection 2.7]{bernig_fu_solanes}) is denoted by $\mathrm{LK}(M)$. Conjecture 3 of \cite{bernig_fu_solanes} states that for every Riemannian manifold $M$ we have 
\begin{displaymath}
 \mathfrak{a}(M) = \mathrm{LK}(M).
\end{displaymath}
Neither of the inclusions seems to be known. 

We study a version of this conjecture on $M=\mathbb{CP}^n_\lambda, \lambda \in \R$. More precisely, we characterize the space of all invariant valuations on $\mathbb{CP}_\lambda^n$ which leave the space of {\it invariant} angular curvature measures invariant.  As in the flat case, the algebra $\mathcal V^n_\lambda$ of $G_\lambda$-invariant valuations in $\mathbb{CP}^n_\lambda$ has two natural generators: $t_\lambda$ and 
$s$ (cf. \cite{bernig_fu_solanes}). 

\begin{Theorem} \label{mainthm_angular}
The algebra of $G_\lambda$-invariant valuations on $\mathbb{CP}_\lambda^n$ leaving the space of angular invariant curvature measures invariant is given by the elements $p(t_\lambda,s)\in\mathcal V^n_\lambda$, with $p\in \mathbb C[t,s]$ such that 
\begin{displaymath}
 tsu \frac{\partial p}{\partial s}\left(\frac{t}{(1-\lambda s)^\frac{1}{2}},s\right)=0 \text{ in } \mathrm{Val}^{\mathrm U(n)}.
\end{displaymath}
\end{Theorem}

If $p$ is a polynomial in $t_\lambda$ alone, then clearly $p$ satisfies the angularity condition. However, since  $\mathrm{Val}^{\mathrm U(n)}$ contains zero divisors, the equation in the corollary does not imply that $\frac{\partial p}{\partial s}=0$ in general. For instance, the image of the polynomial $p(t,s)=t^4-6st^2+6s^2$ is angular in $\Val^{\mathrm U(4)}$, since $tsu(-6t^2+12s)=0$. However, $p$ can not be written as a polynomial in $t$ alone, since there are no relations between $t$ and $s$ of degree $4$.

\section{Background and notations}
\label{sec_background}

For the reader's convenience, we collect here some results from \cite{bernig_fu_solanes} which will be needed in the sequel. 

Let $M$ be a smooth manifold of dimension $n$ which for simplicity is assumed to be oriented and connected. The space of compact submanifolds with corners is denoted by $\mathcal{P}(M)$. A smooth curvature measure is a functional of the form 
\begin{displaymath}
 \Phi(P,\beta):=\int_{N(P) \cap \pi^{-1}\beta} \omega+\int_{P \cap \beta} \rho,
\end{displaymath}
where $N(P)$ is the conormal cycle of $P$ (which is a Legendrian cycle in the cosphere bundle $S^*(M)$); $\beta$ a Borel subset of $M$; $\omega \in \Omega^{n-1}(S^*M)$, $\pi:S^*(M) \to M$ the natural projection and $\rho \in \Omega^n(M)$. The space of smooth curvature measures is denoted by $\mathcal{C}(M)$.

A functional of the form 
\begin{displaymath}
 \phi(P):=\int_{N(P)} \omega+\int_{P} \rho,
\end{displaymath}
is called smooth valuation on $M$. The corresponding space is denoted by $\mathcal{V}(M)$. The obvious map $\glob:\mathcal{C}(M) \to \mathcal{V}(M)$ is called globalization map. 

Alesker \cite{alesker04_product} has introduced a product structure on $\mathcal{V}(M)$, which was generalized in \cite{bernig_fu_solanes} to a module structure of $\mathcal{C}(M)$ over $\mathcal{V}(M)$ satisfying 
\begin{displaymath}
 \glob(\phi\cdot \Phi)=\phi \cdot \glob \Phi, \quad \phi \in \mathcal{V}(M), \Phi \in \mathcal{C}(M).
\end{displaymath}

In the special case where $M=V$ is an affine space, we let $\Curv=\Curv(V)$ be the space of translation invariant smooth curvature measures and $\Val^\infty=\Val^\infty(V)$ the space of translation invariant smooth valuations. Then $\Curv$ is a module over $\Val^\infty$.

Let us now specialize to the case $\Curv^{\mathrm U(n)}$ of unitarily invariant elements in $\Curv$. As a vector space, it is generated by elements $\Delta_{k,q}, 0 \leq k \leq 2n, \max\{0,k-n\} \leq q \leq \left\lfloor\frac{k}{2}\right\rfloor$ and $N_{k,q}, 1 \leq k \leq 2n-3, \max\{0,k-n+1\} \leq q < \frac{k}{2}$. They are defined in terms of invariant differential forms on the sphere bundle, we refer to \cite[Subsection 3.1]{bernig_fu_solanes} for details. The globalization map is injective on $\mathrm{span}\{\Delta_{k,q}\}$ while its kernel is spanned by $\{N_{k,q}\}$.

The algebra $\Val^{\mathrm U(n)}$ is generated by two elements $t,s$. The module structure of $\Curv^{\mathrm U(n)}$ over $\Val^{\mathrm U(n)}$ was computed in \cite{bernig_fu_solanes} as follows. 

\begin{align}
 s \cdot \Delta_{kq} & =\frac{(k-2q+2)(k-2q+1)}{2\pi(k+2)}\Delta_{k+2,q}+\frac{2(q+1)(k-q+1)}{\pi(k+2)}\Delta_{k+2,q+1} \label{eq_mult_s_delta} \nonumber\\
& \quad
-\frac{(k-2q+2)(k-2q+1)}{\pi(k+2)(k+4)}N_{k+2,q}-\frac{2(q+1)(k-2q)}{\pi(k+2)(k+4)}N_{k+2,q+1},\\
 s\cdot N_{kq} & =\frac{(k-2q+2)(k-2q+1)}{2\pi(k+4)}N_{k+2,q}+\frac{2(q+1)(k-q+2)}{\pi(k+4)}N_{k+2,q+1}, \label{eq_mult_s_n}\\
 t\cdot \Delta_{kq} & = \frac{\omega_{k+1}}{\pi \omega_k} \left( (k-2q+1)\Delta_{k+1,q}+2(q+1)\Delta_{k+1,q+1}\right), \label{eq_mult_t_delta}\\
 t\cdot N_{kq} & =\frac{\omega_{k+1}}{\pi \omega_k} \frac{k+2}{k+3} \left((k-2q+1)
N_{k+1,q}+\frac{2(q+1)}{k-2q}(k-2q-1) N_{k+1,q+1}\right) \label{eq_mult_t_n}
\end{align}
 where $\omega_i$ denotes the volume of the $i$-dimensional unit ball.

The map $\phi \mapsto \phi \cdot \Delta_{0,0}$ is called the $\mathfrak l$-map, while the map $\phi \mapsto \phi \cdot N_{1,0}$ is called the $\mathfrak{n}$-map. Then $\Curv^{\mathrm U(n)}$ is generated by the images of the $\mathfrak l$- and $\mathfrak n$-maps. More precisely, 
\begin{equation}\label{curv_decomposition}
 \Curv^{\mathrm U(n)}=\Val^{\mathrm U(n)} \Delta_{0,0} \oplus\widetilde{\Val}^{\mathrm U(n)} N_{1,0},
\end{equation}
where $\widetilde{\Val}^{\mathrm U(n)}:=\C[t,s]/(q_{n-1},q_n)$ is a graded quotient of $\Val^{\mathrm U(n)}$. More precisely, it was shown in \cite{bernig_fu_solanes} that $\widetilde{\Val}^{\mathrm U(n)}:=\C[t,s]/(g_{n-1},g_n)$, where $g_n$ is the degree $n$-part in 
\begin{displaymath}
 e^t \frac{\sin \sqrt{u}-\sqrt{u}\cos \sqrt{u}}{2\sqrt{u}^3}.
\end{displaymath}

By \cite[Lemma 4.4]{wannerer_unitary_module} and \cite[Proposition 5.15]{bernig_fu_solanes}, we have 
\begin{displaymath}
 q_n=\frac{(-1)^{n+1}(n+3)!}{2^n} g_n,
\end{displaymath}
so that the two descriptions of $\widetilde{\Val}^{\mathrm U(n)}$ agree.

The algebra of $G_\lambda$-invariant smooth valuations on $\mathbb{CP}^n_\lambda$ is denoted by $\mathcal{V}^n_\lambda$. It is generated by the generator $t_\lambda$ of the Lipschitz-Killing algebra of $\mathbb{CP}^n_\lambda$ and another element denoted by $s$ \cite{bernig_fu_solanes}.  

By the transfer principle, the space $\mathcal{C}(\mathbb{CP}^n_\lambda)^{G_\lambda}$ of invariant smooth curvature measures on $\mathbb{CP}^n_\lambda$ can be naturally identified with $\Curv^{\mathrm U(n)}$. The local kinematic formulas on the different $\mathbb{CP}^n_\lambda$ are then formally identical. There is a module structure on $\mathcal{C}(\mathbb{CP}^n_\lambda)^{G_\lambda}$ over $\mathcal{V}^n_\lambda$ which is closely related to the so-called semi-local formulas. These are maps 
\begin{align*}
 \bar k_\lambda:\Curv^{\mathrm U(n)} & \to \Curv^{\mathrm U(n)} \otimes \mathcal{V}^n_\lambda\\
 \Phi & \mapsto \left[(P_1,\beta,P_2) \mapsto \int_{G_\lambda} \Phi(P_1 \cap gP_2,\beta)dg\right].
\end{align*}

Writing $\glob_\lambda:\Curv^{\mathrm U(n)} \to \mathcal{V}^n_\lambda$ for the globalization map on $\CP^n_\lambda$, we obviously have $\bar k_\lambda=(\mathrm{id} \otimes \mathrm{glob}_\lambda) \circ K$. The global kinematic formulas on $\mathcal{V}^n_\lambda$ are given by $k_\lambda\circ\glob_\lambda=(\glob_\lambda\otimes\glob_\lambda)\circ K$; i.e.  
\begin{align*}
 k_\lambda:\mathcal{V}^n_\lambda & \to \mathcal{V}^n_\lambda \otimes \mathcal{V}^n_\lambda\\
 \phi & \mapsto \left[(P_1,P_2) \mapsto \int_{G_\lambda} \phi(P_1 \cap gP_2)dg\right]
\end{align*}

Let $\mathrm{pd}_\lambda:\mathcal{V}^n_\lambda \to \mathcal{V}^{n*}_\lambda$ be the normalized Poincar\'e duality, see \cite[Definition 2.14]{bernig_fu_solanes}. We will use the notation $\mathrm{PD}:=\mathrm{pd}_0:\Val^{\mathrm U(n)} \to \Val^{\mathrm U(n)*}$. This map satisfies $\langle \mathrm{PD}\phi,\psi\rangle=(\phi \cdot \psi)_{2n}$, where $\phi \cdot \psi$ is the Alesker product and the subindex $2n$ denotes the component of highest degree.

The module product of $\mathcal{V}^n_\lambda$ on $\Curv^{\mathrm U(n)}$ is a map $\bar m_\lambda \in \mathrm{Hom}(\mathcal{V}^n_\lambda \otimes \Curv^{\mathrm U(n)},\Curv^{\mathrm U(n)})$, which turns out to be  equivalent to the semi-local  kinematic formulas by 
\begin{equation}\label{semi_local_ftaig}
 \bar m_\lambda=(\mathrm{id} \otimes \mathrm{pd}_\lambda) \circ \bar k_\lambda,
\end{equation}
under the identification 
\begin{displaymath}
\mathrm{Hom}(\mathcal{V}^{n}_\lambda \otimes \Curv^{\mathrm U(n)},\Curv^{\mathrm U(n)})= \mathrm{Hom}(\Curv^{\mathrm U(n)},\Curv^{\mathrm U(n)} \otimes \mathcal{V}^{n*}_\lambda). 
\end{displaymath}

Similarly, the Alesker product is a map $m_\lambda:\mathcal{V}^n_\lambda \otimes \mathcal{V}^n_\lambda \to \mathcal{V}^n_\lambda$ which satisfies $(\mathrm{pd}_\lambda \otimes \mathrm{pd}_\lambda) \circ k_\lambda=m_\lambda^* \circ \mathrm{pd}_\lambda$.

We may summarize the different maps in the following commuting diagram. 
\begin{equation} \label{eq_local_and_global_formulas_curved}
\xymatrix{ \Curv^{\mathrm U(n)} \ar[r]^-{K} \ar[d]^{\mathrm{id}} &  \Curv^{\mathrm U(n)} \otimes \Curv^{\mathrm U(n)}  \ar[d]^{\mathrm{id} \otimes \mathrm{glob}_\lambda}\\
\Curv^{\mathrm U(n)} \ar[r]^-{\bar k_\lambda} \ar[d]^{\mathrm{glob}_\lambda} &  \Curv^{\mathrm U(n)} \otimes \mathcal{V}^n_\lambda  \ar[d]^{\mathrm{glob}_\lambda \otimes \mathrm{id}}\\
\mathcal{V}^n_\lambda \ar[r]^-{k_\lambda} \ar[d]^{\mathrm{pd}_\lambda} & \mathcal{V}^n_\lambda \otimes \mathcal{V}^n_\lambda \ar[d]^{\mathrm{pd}_\lambda \otimes \mathrm{pd}_\lambda}\\
\mathcal{V}^{n*}_\lambda \ar[r]^-{m_\lambda^*} & \mathcal{V}^{n*}_\lambda \otimes \mathcal{V}^{n*}_\lambda}
\end{equation}

\section{The algebra structure on $\Curv^{\mathrm{U}(n)*}$}

In this section, we will prove Theorem \ref{mainthm_algebra}. 

Since local and global kinematic formulas are intertwined by  $\glob_\lambda$, the maps $\glob_\lambda^*: \mathcal{V}_\lambda^{n*} \to \Curv^{\mathrm{U}(n)*}$ and $ \glob_\lambda^* \circ \mathrm{pd}_\lambda: \mathcal{V}_\lambda^{n} \to \Curv^{\mathrm{U}(n)*}$ are injective algebra morphisms. Given an element $\phi \in \mathcal{V}_\lambda^{n}$, we will denote by $\bar \phi$ its image in $\Curv^{\mathrm{U}(n)*}$. In particular, the elements $t,s \in \Val^{\mathrm{U}(n)}$ give rise to elements $\bar t,\bar s \in \Curv^{\mathrm{U}(n)*}$.

Let us denote by $\Delta_{k,q}^*,N_{k,q}^* \in \Curv^{\mathrm U(n)*}$ the dual basis of the basis $\Delta_{k,q}, N_{k,q}$ of $\Curv^{\mathrm U(n)}$. The unit element in $\Val$ is the Euler characteristic $\chi$. The unit element in $\Curv^{\mathrm U(n)*}$ is $\bar\chi=\Delta_{2n,n}^*$.

\begin{Proposition} \label{prop_module_vs_dual}
 Let $\phi \in \mathcal{V}_\lambda^{n}$ and $\bar m_\phi:\Curv^{\mathrm U(n)} \to \Curv^{\mathrm U(n)}$ be the module multiplication by $\phi$. Then the dual map
$\bar m_\phi^*:\Curv^{\mathrm U(n)*} \to
\Curv^{\mathrm U(n)*}$ equals multiplication by $\bar \phi$.
\end{Proposition}

\proof
Let $L \in \Curv^{\mathrm U(n)*}$ and $\Phi \in \Curv^{\mathrm U(n)}$. By  \eqref{semi_local_ftaig},
\begin{align*}
\langle \bar m_\phi^*L,\Phi \rangle & =  \langle L, \bar m_\phi(\Phi)\rangle \\
 & =\langle L \otimes \pd_\lambda(\phi),\bar k_\lambda(\Phi)\rangle \\
& = \langle L \otimes \pd_\lambda(\phi),(\mathrm{id} \otimes \glob_\lambda) \circ K(\Phi)\rangle\\
& = \langle L \otimes (\glob_\lambda^* \circ \pd_\lambda)(\phi),K(\Phi)\rangle\\
& = \langle K^*(L \otimes \bar \phi),\Phi\rangle\\
& = \langle  L \cdot \bar \phi,\Phi\rangle.
\end{align*}
\endproof

\begin{Lemma} \label{lemma_dual_product}
 \begin{align*}
  \bar t \Delta_{k,q}^* & = \frac{\omega_k}{\pi \omega_{k-1}} \left[(k-2q)\Delta_{k-1,q}^*+2q\Delta_{k-1,q-1}^*\right], \\
  \bar t N_{k,q}^* & = \frac{\omega_k}{\pi \omega_{k-1}} \frac{k+1}{k+2} \left[(k-2q)N_{k-1,q}^*+\frac{2q(k-2q)}{k-2q+1}N_{k-1,q-1}^*\right],  \\
  \bar s \Delta_{k,q}^* & = \frac{(k-2q)(k-2q-1)}{2\pi k} \Delta_{k-2,q}^*+\frac{2q(k-q)}{\pi k} \Delta_{k-2,q-1}^*,  \\
  \bar s N_{k,q}^* & = -\frac{(k-2q)(k-2q-1)}{\pi k (k+2)}\Delta_{k-2,q}^*-\frac{2q(k-2q)}{\pi k (k+2)} \Delta_{k-2,q-1}^*, \\
  & \quad + \frac{(k-2q)(k-2q-1)}{2\pi (k+2)} N_{k-2,q}^*+\frac{2q(k-q+1)}{\pi(k+2)}N_{k-2,q-1}^*.
\end{align*}
\end{Lemma}
\proof
These equations follow from \eqref{eq_mult_s_delta}, \eqref{eq_mult_s_n}, \eqref{eq_mult_t_delta}, \eqref{eq_mult_t_n} and Proposition \ref{prop_module_vs_dual}.
\endproof

 In particular, $\bar t= \frac{2n\omega_{2n}}{\pi\omega_{2n-1}}\Delta_{2n-1,n-1}^*$ and $\bar s=\frac{n}{\pi}\Delta_{2n-2,n-1}^*$. It will be also useful to compute
\begin{align}\label{eq_stbar}
\bar t \bar s & = \frac{4\omega_{2n}(n-1)n}{\omega_{2n-1}(2n-1)\pi^2}  n
 \Delta_{2n-3,n-2}^*\\
\bar t^3 & = \frac{4\omega_{2n}(n-1)n}{\omega_{2n-1}(2n-1)\pi^2} \left(4 (n-2)
\Delta_{2n-3,n-3}^*+ 6\Delta_{2n-3,n-2}^*\right).\label{eq_t3bar}
 \end{align}

Besides $\bar t,\bar s$ we will need a third element in $\Curv^{\mathrm U(n)*}$, which will be denoted by $\bar v$ (even though it is not the image of an element $v \in \Val^{\mathrm U(n)}$). Namely, 
\begin{equation} \label{eq_def_vbar}
 \bar v := \frac{16\omega_{2n}n(n-1)(n-2)}{\omega_{2n-1}(2n-1)\pi^2} 
\left(\Delta_{2n-3,n-3}^*-\Delta_{2n-3,n-2}^*-\frac{2n-1}{2(n-2)}N_{2n-3,n-2}^*\right).
\end{equation}

As we will see, $\bar v$ vanishes on the image of the $\mathfrak l$-map and is characterized (up to scaling) by this property.

\proof[Proof of Theorem \ref{mainthm_algebra}]
The image of the $\mathfrak{l}$-map in $\Curv^{U(n)}_{2n-3}$ has dimension $2$, while $\dim \Curv^{U(n)}_{2n-3}=3$. We claim that $\bar v$ vanishes on the image
of the $\mathfrak{l}$-map, which defines $\bar v$ uniquely (up to scale). 

By  \cite[Lemma 5.12]{bernig_fu_solanes}, and \eqref{eq_mult_s_delta},  
\begin{align}
\mathfrak{l}(t^{2n-3}) & = \frac{2^{2n-3}(n-2)!}{\pi^{n-1}}(\Delta_{2n-3,n-3}+\Delta_{2n-3,n-2}) \label{lt2n3}\\
\mathfrak{l}(t^{2n-5}u) & =
\frac{2^{2n-3}(n-2)!}{\pi^{n-1}(2n-3)}
\left(\frac{n-3}{n-2}\Delta_{2n-3,n-3}+\Delta_{2n-3,n-2}-\frac{2}{2n-1}N_{2n-3,n-2}\right), \label{lt2n5u}
\end{align}
from which the claim follows. 

Using reverse induction on the degree $k$ and Lemma \ref{lemma_dual_product}, one can show that each element $L \in \Curv^*$ may be written as 
\begin{equation}\label{dual_decomposition}
L=p(\bar s,\bar t) + q(\bar s,\bar t)\bar v
\end{equation}
with polynomials $p,q \in \mathbb{C}[\bar s,\bar t]$. In particular, $\bar s,\bar t,\bar v$ generate
the algebra $\Curv^{\mathrm U(n)*}$.

It follows from \cite{fu06} and \cite[Proposition 5.15]{bernig_fu_solanes} that if $p$ is in the ideal generated by $f_{n+1}(\bar s,\bar t),f_{n+2}(\bar s,\bar t)$ 
and $q$ is in the ideal generated by $q_{n-1}(\bar s,\bar t)$ and $q_n(\bar s,\bar t)$ then $L=p(\bar s,\bar t)+q(\bar s,\bar t)\bar v=0$.  Indeed, for $\Phi=\mu\Delta_{0,0}+\varphi N_{1,0}\in \Curv^{\mathrm{U}(n)}$, 
\[
 \langle L,\Phi\rangle=\langle 1,p(s,t)\Phi\rangle+\langle \bar v,q(s,t)\mu \Delta_{0,0}\rangle +\langle \bar v,q(s,t)\varphi N_{1,0}\rangle=0.
\]

By looking at the dimensions, one sees that there can be no more relations of  degree 1 in $\bar v$. This fixes the algebra
structure on $\Curv^{\mathrm U(n)*}$, except that we have to write $\bar v^2$  in the form \eqref{dual_decomposition}. 

There is an algebra isomorphism $I_\lambda:\Val^{\mathrm U(n)} \to \mathcal{V}^n_\lambda$ \cite[Thm. 3.17]{bernig_fu_solanes} defined by $t \mapsto t_\lambda \sqrt{1-\lambda s}, s \mapsto s$. The map $H_\lambda:=\mathrm{PD}^{-1} \circ I_\lambda^* \circ \mathrm{pd}_\lambda \circ \glob_\lambda:\Curv^{\mathrm U(n)} \to \Val^{\mathrm U(n)}$ was studied in \cite[Section 6]{bernig_fu_solanes}. With $H_0':=\left.\frac{d}{d\lambda}\right|_{\lambda=0} H_\lambda$ we have $H_0' \circ \mathfrak l=D_1, H_0' \circ \mathfrak n=D_2$, where $D_1,D_2:\Val^{\mathrm U(n)} \to \Val^{\mathrm U(n)}$ are covered by the maps 
\begin{displaymath}
 D_1p:=\frac{t^2-2s}{2}p-\frac{tu}{4}\frac{\partial p}{\partial t}, D_2p:=-\frac{3\pi u t}{8}p+\frac{\pi u^2}{8} \frac{\partial p}{\partial t}, p \in \C[t,s].
\end{displaymath}

 Using $$t^{2n-3}s=\frac{n}{2(2n-1)}t^{2n-1},\quad t^{2n-5}s^2=\frac{n(n-1)}{4(2n-1)(2n-3)} t^{2n-1}$$ in $\Val^{\mathrm U(n)}$, one easily computes that
\begin{align}
 \mathfrak{l}(t^{2n-3})+\frac{2(2n-3)}{3\pi} \mathfrak{n}(t^{2n-4}) & \in \ker H_0',\label{first_eq}\\
 \mathfrak{l}(t^{2n-5}u)+\mathfrak{n}\left(\frac{2}{\pi}t^{2n-4}\right) & \in \ker H_0'.\label{second_eq}
\end{align}

By \cite[Lemma 5.12.]{bernig_fu_solanes}, we have 
\begin{equation} \label{nt2n4}
 \mathfrak{n}(t^{2n-4}) = \frac{3}{4} \frac{(2n-4)!\omega_{2n-1}}{\pi^{2n-3}} N_{2n-3,n-2}.
\end{equation}
Using \eqref{lt2n3},\eqref{lt2n5u} and \eqref{nt2n4}, we see that \eqref{first_eq} and \eqref{second_eq} are equivalent to 
\begin{align}
 \Delta_{2n-3,n-3}+\Delta_{2n-3,n-2}+\frac{1}{2n-1} N_{2n-3,n-2} &\in \ker H_0'\\
\frac{n-3}{n-2}\Delta_{2n-3,n-3}+\Delta_{2n-3,n-2}+\frac{1}{2n-1}N_{2n-3,n-2}&
\in \ker H_0'.
\end{align}

It follows that 
\begin{displaymath}
 \Delta_{2n-3,n-2}^*-(2n-1) N_{2n-3,n-2}^* \in \mathrm{Im} (H_0')^*.
\end{displaymath}

This element is (up to a scaling factor) just $\bar v+\bar t \bar u$  (see \eqref{eq_stbar},\eqref{eq_t3bar}). By \cite[Prop. 6.2.]{bernig_fu_solanes}, we have $(H_0' \otimes H_0') \circ K=0$. Dualizing, we obtain that the restriction of the product to $\mathrm{Im}(H_0')^*$ vanishes. In particular, the square of $\bar v+\bar t \bar u$ vanishes. 
\endproof

We may now complete the description of the product structure on $\Curv^{\mathrm{U}(n)*}$ given in Lemma \ref{lemma_dual_product}. 

\begin{Lemma} \label{lemma_dual_product_v}
\begin{align*}
 \bar v \Delta_{k,q}^* & = \frac{16}{\pi^2} \frac{\omega_k}{\omega_{k-1}(k-1)} \bigg[6\binom{q}{3}\Delta^*_{k-3,q-3}+\binom{q}{2}(k-4q+4)\Delta^*_{k-3,q-2}\\
  & \quad -\binom{q}{2}(k-2q) \Delta^*_{k-3,q-1}-\frac{k-1}{k-2q+1}\binom{q}{2}N^*_{k-3,q-2}\bigg],\\
  \bar v N_{k,q}^* & = \frac{\omega_k (k-2q)}{\pi^2 \omega_{k-1} (k+2)} \bigg[\frac{(k-2q-1)(k-2q-2)}{k-1} \Delta_{k-3,q}^*\\
  & \quad +12 \frac{(2k-4q-1)q}{k-1}\Delta^*_{k-3,q-1}+24 \frac{q(q-1)}{k-1} \Delta_{k-3,q-2}^*\\
  & \quad +32 \frac{q-2}{k-2q+3} \binom{q}{2} N_{k-3,q-3}^*\\
  & \quad +16 \frac{k-4q-3}{k-2q+1} \binom{q}{2} N_{k-3,q-2}^*\\
  & \quad -16 \frac{q+2}{q-1} \binom{q}{2} N_{k-3,q-1}^*\bigg].
\end{align*}
\end{Lemma}

\proof
There can be only one linear operator $\sigma$ of degree $-3$ acting on $\Curv^{\mathrm U(n)*}$  with the following properties: $\sigma$ commutes with multiplications by $\bar t$ and $\bar s$, $\sigma \Delta_{2n,n}^*=\bar v$ and 
\begin{equation} \label{eq_char_barv}
 \sigma \bar v+2\bar t \bar u \bar v+(\bar t\bar u)^2=0.
\end{equation}

Indeed, for polynomials $p_1,p_2$ in $\bar t,\bar s$, $\sigma(p_1+p_2\bar v)=p_1 \sigma \Delta_{2n,n}^*+p_2 \sigma \bar v$ is uniquely determined by these properties. The multiplication by $\bar v$ has these properties, where \eqref{eq_char_barv} follows from $(\bar v+\bar t \bar u)^2=0$ which is true by Theorem \ref{mainthm_algebra}. The operator on the right hand side of the displayed equations also has these properties (which is a bit tedious to verify), and hence both sides of the displayed equation agree.
\endproof

\section{The image of $\mathcal{V}^n_\lambda$ in $\Curv^{\mathrm U(n)*}$}

Recall the  injection of algebras $\mathcal{V}^n_\lambda \longrightarrow \Curv^{\mathrm U(n)*}$ given by $\phi\mapsto\bar\phi={\glob_\lambda^*\circ\mathrm{pd}_\lambda(\phi)}$. Since $\mathcal{V}^n_\lambda$ is generated by $s$ and $t_\lambda$, we can describe this morphism by finding the images of $s$ and $t_\lambda$. Since module multiplication by $s$ is independent of the curvature  (see \cite[Prop.5.2]{bernig_fu_solanes}), it follows that $s \in \mathcal{V}^n_\lambda$ is mapped to the same element $\bar s \in \Curv^{\mathrm U(n)*}$ for all $\lambda$. It remains to find the image of $t_\lambda$. 

\begin{Lemma}\label{lemma_pd_ell_en}
For all $p,q \in \C[t,s]$ we have
\begin{equation}\label{eq_pd_ell}
\langle \bar q,\mathfrak l(p)\rangle=\langle \PD(q),p\rangle
\end{equation}  
and 
\begin{equation}\label{eq_pd_en}
 \langle \bar v,\mathfrak n(p)\rangle=\langle \PD(p),e\rangle
\end{equation}
where $e:=-\frac{\pi}{2} u^2$. 
\end{Lemma}

\proof
The first equation follows from 
\begin{align*}
\langle \bar q,\mathfrak l(p)\rangle & = \langle \glob^* \circ \PD(q),\mathfrak l(p)\rangle\\
& = \langle \PD(q),\glob \circ \mathfrak l(p)\rangle\\
& = \langle \PD(q),p\rangle.
\end{align*}

For the second equation, we note that both sides vanish if $p$ is not of degree $2n-4$. We may thus suppose that $p$ is a linear combination of $t^{2n-4}, t^{2n-6}u,t^{2n-8}u^2$. Using \cite[Lemma 5.12]{bernig_fu_solanes} and \cite[Prop. 3.7]{bernig_fu_hig} this is a tedious, but straightforward computation. 
\endproof

\begin{Proposition} \label{prop_image_tlambda}
 The image of $t_\lambda$ in $\Curv^{\mathrm U(n)*}$ is given by 
 \begin{displaymath}
  \bar t_\lambda= \frac{\bar t-\frac{\lambda \bar t^3}{4}}{(1-\lambda \bar s)^\frac{3}{2}}+ \frac{\lambda}{4(1-\lambda \bar s)^\frac{3}{2}} \bar v.
 \end{displaymath}
\end{Proposition}

\proof
Let $p_1:=\frac{\bar t-\frac{\lambda \bar t^3}{4}}{(1-\lambda \bar s)^\frac{3}{2}}, p_2:=\frac{\lambda}{4(1-\lambda \bar s)^\frac{3}{2}}$. We have to show that $\bar t_\lambda=p_1+p_2 \bar v$. 

By \cite[Thm. 6.7]{bernig_fu_solanes} we have 
\begin{align*}
 t_\lambda \Delta_{0,0} & =p_1\Delta_{0,0}+r_1N_{1,0}, \\
 t_\lambda N_{1,0} & =p_2 e \Delta_{0,0}+r_2N_{1,0},
\end{align*}
where $r_1,r_2$ are explicitly known but irrelevant for our purpose.  

Let $q_1 \in \Val^{\mathrm U(n)}$. Then  by \eqref{eq_pd_ell}
\begin{displaymath}
\langle p_1+p_2\bar v,q_1 \Delta_{0,0}\rangle=\langle p_1,q_1 \Delta_{0,0}\rangle= \langle \PD(p_1),q_1\rangle
\end{displaymath}
and 
\begin{displaymath}
 \langle \bar t_\lambda,q_1 \Delta_{0,0}\rangle= \langle\Delta_{2n,n}^*,   t_\lambda q_1 \Delta_{0,0}\rangle=\langle\Delta_{2n,n}^*, q_1(p_1\Delta_{0,0}+r_1N_{1,0})\rangle= \langle \PD(p_1),q_1\rangle.
\end{displaymath}

Next we use  \eqref{eq_pd_en} and compute, for $q_2\in\mathrm{Val}^{\mathrm U(n)}$
\begin{displaymath}
\langle p_1+p_2\bar v,q_2N_{1,0}\rangle=\langle p_2\bar v,q_2N_{1,0}\rangle=\langle \bar v,p_2q_2N_{1,0}\rangle=\langle \PD(p_2q_2),e\rangle
\end{displaymath}
and  
\begin{displaymath}
 \langle t_\lambda,q_2N_{1,0}\rangle=\langle q_2 t_\lambda N_{1,0},\Delta_{2n,n}^*\rangle=\langle q_2(p_2e\Delta_{0,0}+r_2N_{1,0}),\Delta_{2n,n}^*\rangle= \langle \PD(p_2q_2),e\rangle. 
\end{displaymath}
It follows that $\bar t_\lambda$ and $p_1+p_2 \bar v$ act the same on the images of the $\mathfrak l$- and $\mathfrak n$-maps, hence on all $\Curv^{\mathrm U(n)}$, which finishes the proof. 
\endproof

Next we want to describe the image of $\mathcal V_n^\lambda$ in $\Curv^{\mathrm U(n)*}$.

\begin{Lemma} \label{lemma_derivative_t}
 Let $p \in \C[[t,s]]$. 
 \begin{enumerate}
  \item If $p=0$ in $\Val^{\mathrm U(n)}$, then $\frac{\partial p}{\partial t}=0$ in $\widetilde{\Val}^{\mathrm U(n)}$.
  \item If $p=0$ in $\widetilde{\Val}^{\mathrm U(n)}$, then $\frac{\partial (tup)}{\partial t}=0$ in $\widetilde{\Val}^{\mathrm U(n)}$.
 \item The map $p \mapsto Qp:=p+\frac{\lambda}{4(1-\lambda s)}\frac{\partial (tup)}{\partial t}$ is an isomorphism of $\widetilde{\Val}^{\mathrm U(n)}$.
 \end{enumerate}
\end{Lemma}

\proof
A direct computation using \cite[Eq. (38)]{bernig_fu_hig} and \cite[Lemma 4.4]{wannerer_unitary_module}  yields 
\begin{displaymath}
 \frac{\partial f_{n+1}}{\partial t}=-\frac{1}{n+2}(2q_n+tq_{n-1})
\end{displaymath}
for all $n$, from which the first statement follows. The second item follows from  
\begin{displaymath}
 \frac{\partial q_{n-1}}{\partial t} tu=(n+3)t^2q_{n-1}+2t q_n,
\end{displaymath}
which can be shown by induction or by using the defining power series for $q$.

For the third item, we first show that $Qp:=p+\frac{\lambda}{4(1-\lambda s)}\frac{\partial(put)}{t}$ is a surjective map on the space of formal power series in $t,s$. 

Indeed, suppose that $Qp-q=\sum_{k=l}^\infty b_k(s)t^k$ for some $l$ and some formal power series $b_k(s)$. Then
\begin{align*}
Q\left(p-\frac{b_l(s)}{1+(l+1)\frac{\lambda s}{1-\lambda s}}t^l\right)-q \equiv 0 \mod t^l
\end{align*}
Continuing this process with $l$ replaced by $l+1$ and so on, we may construct some power series $p$ with $Qp=q$. 

The space $\widetilde{\Val}^{\mathrm U(n)}$ is a quotient of the space of formal power series in $t,s$, and the map $Q$ induces a map on this quotient by ii). This induced map is still surjective. Since $\widetilde{\Val}^{\mathrm U(n)}$ is finite-dimensional, the map has to be an isomorphism. 

\endproof

\begin{Proposition}
Let $\bar w:=\bar v+\bar t \bar u \in \Curv^{\mathrm U(n)*}$ and $p_1,p_2 \in \C[\bar t,\bar s]$. Then $p_1+p_2 \bar w$ belongs to the image of $\mathcal{V}_\lambda^n$ if and only if 
\begin{equation} \label{eq_image_condition}
\frac{\partial p_1}{\partial t} \cdot \frac{\lambda}{4(1-\lambda s)}=p_2 \quad \text{ in } \widetilde{\Val}^{\mathrm U(n)}.
 \end{equation}
\end{Proposition}

\proof
Suppose that $p_1+p_2\bar w$ belongs to the image of $\mathcal{V}_\lambda^n$. Then there is some polynomial $q$ such that $p_1(\bar t,\bar s)+p_2(\bar t,\bar s)\bar w=q(\bar t_\lambda,\bar s)$. 

By Proposition \ref{prop_image_tlambda}, we have 
\begin{equation}\label{eq_tlambda_w}
 \bar t_\lambda=\frac{\bar t}{\sqrt{1-\lambda \bar s}}+\frac{\lambda}{4(1-\lambda \bar s)^\frac{3}{2}} \bar w.
\end{equation}

Since $\bar w^2=0$, the Taylor expansion of $q$ with respect to $\bar w$ stops at the linear term, i.e. 
\begin{displaymath}
 q(\bar t_\lambda,\bar s)=q\left(\frac{\bar t}{\sqrt{1-\lambda \bar s}},\bar s\right)+q_t\left(\frac{\bar t}{\sqrt{1-\lambda \bar s}},\bar s\right) \frac{\lambda}{4(1-\lambda \bar s)^\frac{3}{2}} \bar w,
\end{displaymath}
where $q_t$ denotes the partial derivative with respect to the first variable. 

Using \eqref{dual_decomposition} and the discussion thereafter, we obtain
\begin{align}
 p_1(t,s)+p_2(t,s)tu & =q\left(\frac{t}{\sqrt{1-\lambda s}},s\right)+q_t\left(\frac{t}{\sqrt{1-\lambda s}},s\right) \frac{\lambda}{4(1-\lambda s)^\frac{3}{2}} tu  \quad \text{ in } \Val^{\mathrm U(n)} \label{eq_image_first_cond}\\
 p_2(t,s) & = q_t\left(\frac{t}{\sqrt{1-\lambda s}},s\right) \frac{\lambda}{4(1-\lambda s)^\frac{3}{2}} \quad \text{ in } \widetilde{\Val}^{\mathrm U(n)} \label{eq_image_second_cond}
\end{align}

Taking the derivative of the first equation with respect to $t$ yields an equation in $\widetilde{\Val}^{\mathrm U(n)}$ by Lemma \ref{lemma_derivative_t}(i). Moreover, applying Lemma \ref{lemma_derivative_t}(ii) to \eqref{eq_image_second_cond} will simplify this equation to 
\begin{displaymath}
\frac{\partial p_1}{\partial t}= q_t\left(\frac{t}{\sqrt{1-\lambda s}},s\right) \frac{1}{\sqrt{1-\lambda s}} \quad \text{ in } \widetilde{\Val}^{\mathrm U(n)}.
\end{displaymath}
Multiplying by $\frac{\lambda}{4(1-\lambda s)}$ and using \eqref{eq_image_second_cond} again then yields \eqref{eq_image_condition}.

We thus obtain that the image of $\mathcal{V}^n_\lambda$ is contained in the space of dual curvature measures satisfying \eqref{eq_image_condition}. 

Let us next compare dimensions. Rewrite $p_1+p_2\bar w=:r_1+r_2\bar v$ with $r_1:=p_1+tup_2,r_2:=p_2$. Then \eqref{eq_image_condition} is equivalent to 
\begin{equation} \label{eq_image_condition_w}
\frac{\partial r_1}{\partial t} \cdot \frac{\lambda}{4(1-\lambda s)}=r_2 + \frac{\lambda}{4(1-\lambda s)} \frac{\partial(r_2 tu)}{\partial t}  \quad \text{ in } \widetilde{\Val}^{\mathrm U(n)}.
\end{equation}

The operator $r_2 \mapsto r_2 + \frac{\lambda}{4(1-\lambda s)} \frac{\partial(r_2 t u)}{\partial t}$ is a bijection on $\widetilde{\Val}^{\mathrm U(n)}$ by Lemma \ref{lemma_derivative_t}(iii). Hence to solve \eqref{eq_image_condition_w}, we can take an arbitrary $r_1 \in \Val^{\mathrm U(n)}$, and $r_2 \in \widetilde{\Val}^{\mathrm U(n)}$ is then uniquely determined. It follows that the dimension of the space of dual curvature measures satisfying \eqref{eq_image_condition} equals $\dim  \Val^{\mathrm U(n)}$, which is the same as the dimension of the image of $\mathcal{V}_\lambda^n$ in $\Curv^{\mathrm U(n)*}$. This shows that these spaces agree. 
\endproof

\section{Angular dual curvature measures}

Let $\Ang^{\mathrm U(n)} \subset \Curv^{\mathrm U(n)}$ be the subspace of angular curvature measures  (see \cite[Definition 2.26]{bernig_fu_solanes}). By \cite[Proposition 3.2]{bernig_fu_solanes}, this is the space generated by the curvature measures $\Delta_{k,q}$. 

\begin{Definition}
 A dual curvature measure $L \in \Curv^{\mathrm U(n)*}$ is called angular, if $L \Ang^{\mathrm U(n)\perp} \subset \Ang^{\mathrm U(n)\perp}$, where
 \begin{displaymath}
  \Ang^{\mathrm U(n)\perp}:=\{\Psi \in \Curv^{\mathrm U(n)*}| \langle \Psi,C\rangle=0, \forall C \in \Ang^{\mathrm U(n)}\}
 \end{displaymath}
 is the annihilator of $\Ang^{\mathrm U(n)}$.
\end{Definition}

\begin{Lemma} \label{lemma_partial_int}
Let $p \in \C[t,s], g \in \C[t,u]$ and let $\PD$ be the Alesker-Poincar\'e  duality in $\Val^{\mathrm U(n)}$. Then
 \begin{displaymath}
\left\langle \PD\left(\frac{\partial g}{\partial u}\right),p su \right\rangle =\frac18 \left\langle  \PD(g),(-t^2+(2n-1)u)p-2 \frac{\partial p}{\partial s} us \right\rangle.      
\end{displaymath}
\end{Lemma}

\proof
It suffices to check the formula for polynomials of the form $g=t^{k-2r}s^r$ (written in terms of $t,u$) and $p=t^{l-2i}s^i$ with $k+l=2n-2$. Using  \cite[Eq. (2.57)]{fu_barcelona},
the left hand side equals $ \omega_{2n}^{-1}\left(r\binom{2n-2(r+i+1)}{n-(r+i+1)}-\frac{r}{4} \binom{2n-2(r+i))}{n-(r+i)}\right)$, while the right hand side is $ \omega_{2n}^{-1}\left(\frac{-2n+2i}{8} \binom{2n-2(r+i)}{n-(r+i)}+\frac{2n-2i-1}{2} \binom{2n-2(r+i+1)}{n-(r+i+1)}\right)$. These two expressions agree. 
\endproof

\begin{Lemma} \label{lemma_v_nstar}
Set $\bar r_n:=\frac{4n-2}{n}\bar t \bar s-\bar t^3$. Then $\bar v+\bar r_n$ is a non-zero multiple of $N^*_{2n-3,n-2}$.
\end{Lemma}

\proof
Follows from the explicit formulas \eqref{eq_stbar},\eqref{eq_t3bar}, and \eqref{eq_def_vbar} . 
\endproof

\begin{Lemma} \label{lemma_restriction_n2n3}
 A dual curvature measure $L$ is angular if and only if $L N^*_{2n-3,n-2} \in \Ang^\perp$. 
\end{Lemma}

\proof
Since the space $\Ang^{\mathrm U(n)}$ is spanned by the $\Delta_{k,q}$, its annihilator $\Ang^{\mathrm U(n) \perp}$ is spanned by the $N_{k,q}^*$. Hence $L$ is angular if and only if $LN_{k,q}^*$ is a linear combination of $N_{l,i}^*$'s. This shows the ``only if'' part. 

By \cite[Theorem 6.9]{bernig_fu_solanes} we know that $\bar t$ and $\bar t_\lambda$ are angular. Hence $\bar t_0':=\left.\frac{d}{d\lambda}\right|_{\lambda=0}\bar t_\lambda=-\frac14 \bar t^3+\frac32 \bar t\bar s +\frac14 \bar v$ is angular. 

By reverse induction on $l$, one can show, using Lemmas \ref{lemma_dual_product} and \ref{lemma_dual_product_v}, that every element $N_{l,i}^*$ can be written as $q( \bar t,\bar t_0')N_{2n-3,n-2}^*$ for some polynomial $q$. Now suppose that $LN_{2n-3,n-2}^* \in \Ang^\perp$. Then $LN_{l,i}^*=Lq( \bar t,\bar t_0')N_{2n-3,n-2}^*=q( \bar t,\bar t_0') LN_{2n-3,n-2}^*$ is in $\Ang^\perp$ for every $l,i$, as claimed in the ``if'' part. 
\endproof

The next theorem describes the space $\Ang^{\mathrm U(n)\perp}$.

\begin{Theorem} \label{thm_angular_dual_measures}
 Let $p_1,p_2$ be polynomials in $t, s$. Then $\bar p_1+\bar p_2\bar v$ is angular if and only if 
\begin{displaymath}
tsu\left(-6tp_2+\frac{\partial p_1}{\partial s}-\frac{\partial p_2}{\partial s} tu\right)=0 \quad \text{ in } \Val^{\mathrm U(n)}.
\end{displaymath}
Equivalently, with $\bar w:=\bar v+\bar t \bar u$, $\bar p_1+\bar p_2\bar w$ is angular if and only if
\begin{displaymath}
 tsu\left(\frac{\partial p_1}{\partial s}-2tp_2\right)=0 \text{ in } \Val^{\mathrm U(n)}.
\end{displaymath}

\end{Theorem}

\proof
Let $q_1,q_2$ be polynomials in $t,s$. Then, by Lemma \ref{lemma_v_nstar}, Theorem \ref{mainthm_algebra}, and Lemma \ref{lemma_pd_ell_en} (and $\cong$ meaning equal up to some non-zero factor)
\begin{align*}
 \langle (\bar p_1+\bar p_2 \bar v) &  N^*_{2n-3,n-2}, \mathfrak l(q_1)+\mathfrak n(q_2)\rangle \cong \left\langle (p_1+p_2 \bar v) (\bar v+r_n), \mathfrak l(q_1)+\mathfrak n(q_2)\right\rangle \\
 & = \left\langle \bar p_1 \bar r_n-\bar p_2\bar t^2 \bar u^2+(\bar p_1+\bar p_2 \bar r_n-2 \bar p_2 \bar t \bar u)\bar v, \mathfrak l(q_1)+\mathfrak n(q_2)\right\rangle\\
& = \langle \PD(\tilde p_1),q_1\rangle+ \langle \PD(\tilde p_2e),q_2\rangle,
\end{align*}
where $\tilde p_1=p_1r_n+p_2t^2u^2$ and $\tilde p_2=p_1+p_2r_n-2p_stu$.

By Lemma \ref{lemma_restriction_n2n3}, $\bar p_1+\bar p_2 \bar v$ is angular if and only if $(\bar p_1+\bar p_2 \bar v)N^*_{2n-3,n-2} \in \Ang^{\mathrm U(n)\perp}$. By \cite[Prop. 6.8.]{bernig_fu_solanes}, the curvature measure $\mathfrak l(q_1)+\mathfrak n(q_2)$ is angular if and only if we can write 
\begin{displaymath}
 q_1=g+2u \frac{\partial g}{\partial u}, \quad q_2=\frac{4t}{\pi} \frac{\partial g}{\partial u}
\end{displaymath}
with $g \in \R[t,u]$.
 
It follows that $\bar p_1+\bar p_2\bar v$ is angular if and only if 
\begin{align*}
  0 & =\left\langle \PD(\tilde p_1),g+2u \frac{\partial g}{\partial u}\right\rangle+\left\langle \PD(e \tilde p_2) ,\frac{4t}{\pi} \frac{\partial g}{\partial u}\right\rangle\\
& = \left\langle \PD(\tilde p_1),g\right\rangle+\left\langle \PD(2u\tilde p_1+\frac{4t}{\pi}e \tilde p_2),\frac{\partial g}{\partial u}\right\rangle\\
& = \left\langle \PD(\tilde p_1),g\right\rangle+\left\langle \PD(2u(\tilde p_1-tu \tilde p_2)),\frac{\partial g}{\partial u}\right\rangle
\end{align*}
for each polynomial $g \in \R[t,u]$.

Plugging in the values for $r_n$ and $\tilde p_1,\tilde p_2$, we find that 
\begin{displaymath}
\tilde p_1-tu \tilde p_2=(p_1-p_2tu)(r_n-tu)=-\frac{2}{n} (p_1-p_2tu)ts.
\end{displaymath}

Using Lemma \ref{lemma_partial_int} we obtain the necessary and sufficient condition for angular dual measures:
\begin{multline*}
0=\langle \PD(\tilde p_1),g\rangle-\frac{1}{2n} \left\langle \PD(t(-t^2+(2n-1)u)(p_1-p_2tu)-2t\frac{\partial(p_1-p_2tu)}{\partial s} us),g\right\rangle \quad \\ \forall g \in \R[t,u]. 
\end{multline*}

By the  injectivity of the Alesker-Poincar\'e duality, this is equivalent to 
\begin{align*}
0 & = \tilde p_1-\frac{1}{2n}(t(-t^2+(2n-1)u)(p_1-p_2tu)-2t\frac{\partial(p_1-p_2tu)}{\partial s} us) \\
& = \frac{tsu}{n}\left(-6tp_2+\frac{\partial p_1}{\partial s}-\frac{\partial p_2}{\partial s} tu\right) 
\end{align*}
in $\Val^{\mathrm U(n)}$.
\endproof

\proof[Proof of Theorem \ref{mainthm_angular}]
 Recalling \eqref{eq_tlambda_w}, and $\bar w^2=(\bar v+\bar t\bar u)^2=0$, we obtain that 
\begin{align*}
 p(\bar t_\lambda,\bar s)=p\left(\frac{\bar t}{(1-\lambda \bar s)^\frac{1}{2}},\bar s\right)+\frac{\partial p}{\partial t} \left(\frac{\bar t}{(1-\lambda \bar s)^\frac{1}{2}},\bar s\right) \frac{\lambda}{4(1-\lambda \bar s)^\frac{3}{2}}\bar w.
\end{align*}
The statement then follows from Theorem \ref{thm_angular_dual_measures}.
\endproof

\def\cprime{$'$}

%\bibliographystyle{plain}
%\bibliography{../biblio}
\end{document}